\documentclass[12pt]{amsart}
\usepackage{amsmath}
\usepackage{amsfonts}
\usepackage{amssymb}

\usepackage{verbatim}

\newtheorem{theorem}{Theorem}
\newtheorem{proposition}[theorem]{Proposition}

\newtheorem{lemma}[theorem]{Lemma}

\theoremstyle{definition}
\newtheorem{definition}[theorem]{Definition}
\newtheorem{example}[theorem]{Example}
\newtheorem*{remark}{Remark}

\newcommand{\mcA}{\mathcal{A}}

\newcommand{\mcH}{\mathcal{H}}

\newcommand{\mcJ}{\mathcal{J}}
\newcommand{\mcR}{\mathcal{R}}
\newcommand{\mcS}{\mathcal{S}}
\newcommand{\bbN}{\mathbb{N}}
\newcommand{\bbR}{\mathbb{R}}
\newcommand{\bbZ}{\mathbb{Z}}
\newcommand{\BH}{\mathcal{B}(\mathcal{H})}
\newcommand{\rad}{\mathrm{Rad}}
\newcommand{\prad}{\mathrm{PRad}}

\newcommand{\cstar}{C$^*$}

\newcommand{\bo}[1]{\mathbf{#1}}

\newcommand{\bom}{\mathbf{m}}
\newcommand{\bon}{\mathbf{n}}

\newcommand{\gS}{{\Sigma}}

\begin{document}
\title[The Jacobson radical for analytic crossed products]{The
    Jacobson radical for \\ analytic crossed products}
\author[Donsig]{Allan P. Donsig}
\address{Mathematics \& Statistics Dept.\\
     University of Nebraska---Lincoln\\
     Lincoln, NE,  68588, U.S.A.}
\email{adonsig@math.unl.edu}
\author[Katavolos]{Aristides Katavolos}
\address{Dept.\ of Mathematics\\
    University of Athens\\
    Panepistimioupolis\\
    GR-157 84, Athens\\
    GREECE}
\email{akatavol@eudoxos.math.uoa.gr}
\author[Manoussos]{Antonios Manoussos}
\address{123, Sapfous St.\\
        176 75 Kallithea\\
        Athens\\
    GREECE}
\email{amanou@cc.uoa.gr}
\thanks{To appear, \textit{Journal of Functional Analysis.}}
\thanks{2000 {\it Math.\ Subject Classification.} Primary: 47L65;
    Secondary: 46H10, 46K50.}

\begin{abstract}
We characterise the Jacobson radical of an analytic crossed
product $C_0(X) \times_\phi \bbZ_+$, answering a question first
raised by Arveson and Josephson in 1969. In fact, we characterise
the Jacobson radical of analytic crossed products $C_0(X)
\times_\phi \bbZ^d_+$. This consists of all elements whose
`Fourier coefficients' vanish on the recurrent points of the
dynamical system (and the first one is zero). The
multi-dimensional version requires a variation of the notion of
recurrence, taking into account the various degrees of freedom.
\end{abstract}

\maketitle

There is a rich interplay between operator algebras
and dynamical systems, going back to the founding work of Murray
and von Neumann in the 1930's. Crossed product constructions
continue to provide fundamental examples of von Neumann algebras
and \cstar-algebras. Comparatively recently, Arveson~\cite{Arv67a}
in 1967 introduced a nonselfadjoint crossed product construction,
called the analytic crossed product or the semi-crossed product,
which has the remarkable property of capturing all of the
information about the dynamical system.

The construction starts with a dynamical system $(X,\phi)$,
i.e., a locally compact Hausdorff space $X$ and a continuous, proper
surjection $\phi : X \to X$.
Regarding the elements of $\ell^1(\bbZ_+,C_0(X))$ as formal series
$\sum_{n \ge 0} U^n f_n$, define a multiplication by requiring
$f U = U (f\circ\phi)$.
The analytic crossed product, $C_0(X) \times_\phi \bbZ_+$,
is a suitable completion of $\ell^1(\bbZ_+,C_0(X))$; we give a
detailed discussion below.
Then the property mentioned above is that, subject to a mild
condition on periodic points, two analytic crossed product
algebras are isomorphic as complex algebras if, and only if,
the underlying dynamical systems are topologically conjugate,
i.e., there is a homeomorphism between the spaces which intertwines
the two actions.
In this generality, the result is due to Hadwin and
Hoover~\cite{HadHoo88,HadHoo89}--see also~\cite{Pow92b} which
gives an elegant direct proof of this if the maps $\phi$ are
homeomorphisms and extends the result to analytic crossed products
by finitely many distinct commuting homeomorphisms on $X$,
i.e., by $\bbZ_+^d$.

Arveson's original work~\cite{Arv67a} was for weakly-closed operator
algebras and Arveson and Josephson in~\cite{ArvJos69} gave an
extension to norm-closed operator algebras, including a structure
theorem for bounded isomorphisms between two such algebras.
Motivated by this, they asked if the analytic crossed product algebras
were always semisimple (which would imply that all isomorphisms
are bounded), noting that the evidence suggested a negative answer.
This question stimulated considerable work on the ideal structure of
analytic crossed products.

Another stimulus is the close connections between the ideal
structure of \cstar-crossed products and dynamical systems, such
as the characterisation of primitive ideals of \cstar-crossed
products in terms of orbit closures by Effros and
Hahn~\cite{EffrosHah67}. In this connection, we should mention
Lamoureux's development of a generalisation of the primitive ideal
space for various nonselfadjoint operator algebras, including
analytic crossed products~\cite{Lam93,Lam96a}.

We state our main result for the case $d=1$. Recall that a point
$x \in X$ is \emph{recurrent} for the dynamical system $(X,\phi)$
if for every neighbourhood $V$ of $x$, there is $n \ge 1$ so that
$\phi^n(x) \in V$. If $X$ is a metric space, then this is
equivalent to having a sequence $(n_k)$ tending to
infinity so that $\phi^{n_k}(x)$ converges to $x$.
Let $X_r$ denote the recurrent points of $(X,\phi)$. Denoting
elements of the analytic crossed product by formal series
$\sum_{n\ge 0} U^n f_n$ we prove:
\begin{theorem}
\label{fulrad} If $X$ is a locally compact metrisable space, then
\begin{align*}
\rad( C_0(X)\times_\phi\bbZ_+)
    &= \left\{ \sum_{n \ge 1} U^n f_n \in C_0(X)\times_\phi\bbZ_+:
     f_n |_{X_r}=0\text{ for all }n \right\}. \\
\end{align*}
\end{theorem}

Important progress towards a characterisation has been made by a
number of authors. In~\cite{Muh83}, Muhly gave two sufficient
conditions, one for an analytic crossed product to be semisimple
and another for the Jacobson radical to be nonzero. The sufficient
condition for a nonzero Jacobson radical is that the dynamical
system $(X,\phi)$ possess a \textit{wandering set}, i.e., an open
set $V \subset X$ so that $V,\phi^{-1}(V),\phi^{-2}(V),\ldots$ are
pairwise disjoint. If there are no wandering open sets, then the
recurrent points are dense, so it turns out that this sufficient
condition is also necessary.

Peters in~\cite{Pet84,Pet88a} characterised the strong radical
(namely, the intersection of the maximal (modular, two-sided)
ideals) and the closure of the prime radical and
described much of the ideal structure for analytic crossed
products arising from \emph{free} actions of $\bbZ^+$. He also
gave a sufficient condition for semisimplicity and showed that
this condition is necessary and sufficient for semisimplicity of
the norm dense subalgebra of polynomials in the analytic crossed
product.

Most recently, Mastrangelo, Muhly and Solel~\cite{MasMuhSol94},
using powerful coordinate methods and the crucial idea
from~\cite{Don93b}, characterised the Jacobson radical for
analytic subalgebras of groupoid \cstar-algebras. For those
analytic crossed products which can be coordinatised in this way
(those with a free action), their characterisation is the same as
ours. The asymptotic centre of the dynamical system which is used
in~\cite{MasMuhSol94} is also important to our approach.

However, we are able to dispose of the assumption of freeness (and
thus our dynamical systems can have fixed points or periodic
points); in fact, our methods are applicable to irreversible
dynamical systems having several degrees of freedom (that is, actions
of $\bbZ_+^d$). In the multi-dimensional case the usual notions of
recurrence and centre are not sufficient to describe the Jacobson
radical, as we show by an example. Accordingly, we introduce
appropriate modifications.

After discussing the basic properties of analytic crossed products
and some of the radicals of Banach algebras, we develop the key
lemma in Section~1.
This lemma, which is based on the idea
of~\cite[Lemma~1]{Don93b}, relates (multi-) recurrent points in the
dynamical system with elements not in the Jacobson radical.
In Section~2, we give a characterisation of semisimplicity.
The proof has three ingredients: the key lemma, a sufficient condition
for an element to belong to the prime radical (a descendant of Muhly's
condition mentioned earlier), and a basic fact from dynamical systems
theory which is known in the one-dimensional case.
Our main result, Theorem~\ref{fulrad2}, is proved in the last section
using a modification of the centre of a dynamical system.

\subsection{Definition of Analytic Crossed Products}
Analytic crossed products or semicrossed products have been defined
in various degrees of generality by several authors (see for example
\cite{HadHoo88,Lam96a,Pet84,Pet88a,Pow92b}), generalising the concept
of the crossed product of a \cstar-algebra by a group of *-automorphisms.
To fix our conventions, we present the definition in the form that
we will use it.
Let $X$ be a locally compact Hausdorff space and
$\Phi = \{ \phi_\bon : \bon = (n_1,n_2, \ldots, n_d) \in \bbZ_+^d \}$
be a semigroup of continuous and proper surjections isomorphic
(as a semigroup) to $\bbZ_+^d$.

An action of $\gS =\bbZ_+^d$ on $C_0(X)$ by
isometric *-endomorphisms $\alpha_\bon \, (\bon \in \gS)$ is
obtained by defining $\alpha_\bon(f)=f\circ\phi_\bon$.

We write elements of the Banach space $\ell^{1}(\gS,C_0(X))$ as formal
multi-series $A=\sum_{\bon \in \gS } U_\bon f_\bon$ with the norm
given by $\| A \|_1=\sum \| f_\bon \|_{C_0(X)}$. The multiplication
on $\ell^1(\gS,C_0(X))$ is defined by setting
\[ U_\bon f\, U_\bom g = U_{\bon+\bom} (\alpha_\bom(f)g) \]
and extending by linearity and continuity. With this
multiplication, $\ell^1 (\gS,C_0(X))$ is a Banach algebra.

We will represent $\ell^1 (\gS,C_0(X))$ faithfully as a (concrete)
operator algebra on Hilbert space, and define the analytic crossed
product as the closure of the image.

Assuming we have a faithful action of $C_0(X)$ on a Hilbert space
$H_o$, we can define a faithful contractive representation $\pi$
of $\ell^1(\gS,C_0(X))$ on the Hilbert space
$\mcH=H_{o}\otimes\ell^{2}(\gS)$
by defining $\pi(U_{\bon}f)$ as follows:
\[ \pi(U_{\bon} f) (\xi\otimes e_{\bo{k}}) =
    \alpha_{\bo{k}}(f)\xi\otimes e_{\bo{k+n}} \]

To show that $\pi$ is faithful, let
$A=\sum_{\bon\in\gS} U_\bon f_\bon$ be in $\ell^1(\bbZ_+^d,C_0(X))$
and $x,y\in H_{o}$ be unit vectors.
Since $\pi$ is clearly contractive, the series
$\pi(A)=\sum_{\bon\in\gS} \pi(U_\bon f_\bon)$ converges absolutely.
For $\bom \in \gS$, we have
\begin{align*}
\langle\pi(A)(x\otimes e_{\bo{0}}),y\otimes e_{\bom}\rangle
 &=\sum_{\bon}\langle\pi(U_{\bon} f_{\bon})(x\otimes e_{\bo{0}}),
    y\otimes e_{\bom}\rangle\\
 &=\sum_{\bon} \langle f_\bon x
    \otimes e_\bon, y\otimes e_{\bom} \rangle \\
 &=\langle f_{\bom} x\otimes e_{\bom},
    y\otimes e_\bom\rangle = \langle f_{\bom} x,y\rangle
\end{align*}
as $x\otimes e_\bon$ and $y\otimes e_\bom$ are orthogonal
for $\bon\neq\bom$.
It follows that
\[
\Vert\pi(A)\Vert\geq\Vert f_{\bom}\Vert.
\]
Hence if $\pi(A)=0$ then $f_{\bom}=0$ for all $\bom$, showing $A=0$.
Thus $\pi$ is a monomorphism.

\begin{definition}
The \emph{analytic crossed product}
$\mcA = C_0(X)\times_\phi\bbZ_+^d$ is the closure of the image of
$\ell^1(\bbZ_+^d,C_0(X))$ in $\BH$ in the representation
just defined.
\end{definition}

This is a generalisation of the definition given in~\cite{Pet88a}.
Note that $\mcA$ is in fact independent of the faithful action of
$C_0(X)$ on $H_o$ (up to isometric isomorphism).
\medskip

For $A=\sum U_{\bon}f_{\bon} \in\ell^1(\gS,C_0(X))$ we call
$f_{\bon}\equiv E_{\bon}(A)$ the \emph{${\bon}$-th Fourier
coefficient }of $A$. We have shown above that the maps
$E_{\bon}:\ell^1(\gS,C_0(X))\rightarrow C_0(X)$ are contractive
in the (operator) norm of $\mcA$, hence they extend to
contractions $E_{\bon}:\mcA\to C_0(X)$.

Moreover,
\[
U_{\bom} E_{\bom}(A) =
\frac{1}{(2 \pi)^d} \int_{([- \pi ,
\pi])^d}\theta_{\bo{t}}(A) \exp (-i\bo{m.t})d\bo{t}
\]
where $\bo{m.t} = m_1t_1+ \ldots + m_dt_d$ and the automorphism
$\theta_{\bo{t}}$ is defined first on the dense subalgebra
$\ell^1(\gS,C_0(X))$ by
\[
\theta_{\bo{t}} \left(\sum U_{\bon}f_{\bon} \right)
    = \sum U_{\bon} (\exp (i \bo{t.n}) f_{\bon})
\]
and then extended to $\mcA$ by continuity.

Thus, by injectivity of the Fourier transform on $C(([-\pi,\pi])^d)$,
if a continuous linear form $\eta$ on $\mcA$ satisfies
$\eta (E_{\bom}(A))=0$ for all $\bom$ then (the function
$\bo{t} \to \eta( \theta_{\bo{t}}(A))$ vanishes and hence)
$\eta(A)=0$.
The Hahn-Banach Theorem yields the following remark.

\begin{remark}\label{fej}
Any $A \in \mcA$ belongs to the closed linear span of the set
$\{U_{\bom} E_{\bom}(A) : \bom \in \gS \}$ of its
`associated monomials'.

In particular, $\mcA$ is the closure of the subalgebra $\mcA_0$ of
trigonometric polynomials, i.e., finite sums of monomials.
\end{remark}

As $\theta_{\bo{t}}$ is an automorphism of $\mcA$, we conclude
that if $\mcJ\subseteq \mcA$ is a closed automorphism invariant
ideal (in particular, the Jacobson radical) then for all $B\in\mcJ$
and $\bom\in\gS$ we obtain $U_{\bom}E_{\bom}(B) \in\mcJ$.
Thus, an element $\sum U_{\bon}f_{\bon} $ is in $\mcJ$ if and only
if each monomial $U_{\bon}f_{\bon} $ is in $\mcJ$;
this was first observed (for $d=1$) in~\cite[Proposition~2.1]{Muh83}.
It now follows from the Remark that any such ideal is the closure
of the trigonometric polynomials it contains.

\subsection{Radicals in Banach Algebras}
Recall that an ideal $\mcJ$ of an algebra $\mcA$ is said to be
\emph{primitive} if it is the kernel of an (algebraically) irreducible
representation.
The intersection of all primitive ideals of $\mcA$ is the
\emph{Jacobson radical} of $\mcA$, denoted $\rad\mcA$.

An ideal $\mcJ$ is \emph{prime} if it cannot factor as the product
of two distinct ideals, i.e., if $\mcJ_1,\mcJ_2$ are ideals of
$\mcA$ such that $\mcJ_1\mcJ_2\subseteq\mcJ$ then either
$\mcJ_1\subseteq\mcJ$ or $\mcJ_2\subseteq\mcJ$. The intersection
of all prime ideals is the \emph{prime radical} of $\mcA$, denoted
$\prad\mcA$. An algebra $\mcA$ is \emph{semisimple} if
$\rad\mcA=\{0\}$ and \emph{semiprime} if $\prad\mcA=\{0\}$,
or equivalently, if there are no (nonzero) nilpotent
ideals.

As a primitive ideal is prime, $\prad\mcA \subseteq \rad\mcA$.
Thus a semisimple algebra is semiprime. If $\mcA$ is a
Banach algebra, then the Jacobson radical is closed; indeed every
primitive ideal is the kernel of some \emph{continuous}
representation of $\mcA$ on a Banach space. In fact an
element $A \in\mcA$ is in $\rad\mcA$ if and only if the spectral
radius of $AB$ vanishes for all $B \in\mcA$.

The prime radical need not be closed; it is closed if and only if
it is a nilpotent ideal
(see~\cite{Gra69} or \cite[Theorem~4.4.11]{Palmer94}).
Thus for a general Banach algebra,
$ \prad\mcA \subseteq \overline{\prad\mcA} \subseteq \rad\mcA$.

\section{Recurrence and Monomials}

Our main results will be proved for metrisable dynamical systems;
hence we make the blanket assumption that \emph{$X$ will be a
locally compact metrisable space}. As in the one-dimensional case,
we say that a point $x \in X$ is \emph{recurrent} for the
dynamical system $(X,\Phi)$ if there exists a sequence $(\bon_k)$
tending to infinity so that $\phi_{\bon_k}(x)\to x$.
We will need the following variant:

\begin{definition}
Let $J\subseteq\{1,2,\ldots d\}$. Say $x\in X$ is
$J$-\emph{recurrent} if there exists a sequence $(\bon_{k})$
which is strictly increasing in the directions of $J$ $($that is,
the $j$-th entry of $\bon_{k+1}$ is greater than the $j$-th entry of
$\bon_{k}$ for every $j \in J$ and $k\in\bbN)$ such
that $\lim_{k}\phi_{\bon_{k}}(x)=x$.
Denote the set of all $J$-recurrent points by $X_{Jr}$.

We say that a point $x\in X$ is \emph{strongly recurrent}
if it is $\{1,2,\ldots,d\}$-recurrent.
Finally, $\Sigma_J$ denotes
$\{ \bo{n} \in \bbZ_+^d : \bo{n}_j > 0 \hbox{ for all } j\in J\}$.
\end{definition}

In the multi-dimensional case, the Jacobson radical cannot be
characterised in terms of either the recurrent points (in the
traditional sense) or the strongly recurrent points. To justify
this, we give the following example.

\begin{example}
Let $X=X_0 \cup X_1 \cup X_2$ where $X_i = \bbR \times \{i\}$.
Consider the dynamical system $(X,(\phi_1,\phi_2))$,
where $\phi_1$ acts as translation by $1$ on $X_1$ and as
the identity on $X_0 \cup X_2$ while $\phi_2$ acts as translation
by $1$ on $X_2$ and as the identity on $X_0 \cup X_1$.
It is easy to see that the set of $\{1\}$-recurrent points is $X_0 \cup X_2$,
the set of $\{2\}$-recurrent points is $X_0 \cup X_1$ and the set of
strongly recurrent points is $X_0$.

Choose small neighbourhoods $V_1 \subseteq X_1$ and $V_2 \subseteq X_2$
of $(0,1)$ and $(0,2)$ respectively such that
$\phi_1(V_1)\cap V_1=\emptyset$ and $\phi_2(V_2)\cap V_2=\emptyset$.
Let $f\in C_{0}(X)$ be any function supported on
$V_1\cup V_2$ such that $f(0,1)=f(0,2)=1$.

Then one can verify (as in the proof of Lemma~\ref{prad} in the
next section) that $U_1U_2 f$ is in the prime radical. On the
other hand, neither $U_1 f$ nor $U_2 f$ belong to the Jacobson
radical (they are not even quasinilpotent).

Here, the associated semicrossed product has nonzero Jacobson
radical, although every point is recurrent. Also, the monomial
$U_1 f$ is not in the Jacobson radical, although $f$ vanishes on
the strongly recurrent points. The next lemma shows that for such
a monomial to be in the Jacobson radical, $f$ must vanish on the
$\{1\}$-recurrent points.
\end{example}

The main result of this section is the following lemma,
which is crucial to our analysis.

\begin{lemma}
\label{rad} Let $U_{\bo{q}}f \in\rad(C_0(X)\times_\phi
\bbZ_+^d)$. If $J$ contains the support of $\bo{q}$, then $f$
vanishes on each $J$-recurrent point of $(X,\Phi)$.
\end{lemma}

In order to prove this lemma, we need a basic property of
recurrent points, adapted to our circumstances.
\medskip

\begin{definition}
Given a sequence $\bar{\bon}=(\bon_{k})\subseteq\bbZ_+^d$
we define recursively the \emph{family of indices associated to}
$\bar{\bon}$, denoted $\mcS(\bar{\bon})=(S_0,S_1,S_2,\ldots)$.
as follows: $S_0=\{\bo{0}\},S_1=\{\bon_1\}$ and generally
\[
S_{k+1}=\{\bon_{k+1}+\bom_k+\bo{j}:
    \bo{j}\in\cup_{i=0}^{k}S_i\}
\]
where $\bom_0=\bo{0}$ and
$\bom_k=\bon_k+2\bom_{k-1}$.
\end{definition}

The sets in $\mcS(\bar{\bon})$ will be needed in the
proof of Lemma \ref{rad}: they are the indices of
$\phi$ occurring in the simplification of the inductive
sequence of products given by $P_1=U_{\bon_1} g$ and
$P_{k} = P_{k-1}(U_{\bon_{k}} (g/2^{k-1}) ) P_{k-1}$.
We should also point out that $\cup_i S_i$ is
an IP-set (see~\cite[Section~8.4]{Furstenberg81}) and the next
lemma is a variant on~\cite[Theorem~2.17]{Furstenberg81}.

Recall $\Sigma_J$ denotes $\{(n_1,n_2,\ldots,n_{d}): n_{j}\neq0 \;
\mbox{ for all } j \in J \}$.
Let $\Delta_J$ be the subset of $\Sigma_J$ with entries
in the directions of $J^c$ identically zero.

\begin{lemma} \label{rec}
Let $x$ be in $X$, $J$ be a subset of $\{1,2,\ldots,d\}$.
Suppose that $\lim_{k}\phi_{\bo{p}_{k}}(x)=x$, where $(\bo{p}_{k})$
is a sequence whose restriction to $J$ is strictly increasing while its
restriction to $J^{c}$ is constant.

For each open neighbourhood $V$ of $x$ and each $k\in \mathbb{N}$,
there is $\bo{n}_{k}\in \Delta _{J}$ and $x_{k}\in V$ with
\[
\phi _{\mathbf{s}}(x_{k})\in V\quad \hbox{ for all }\quad
\mathbf{s}\in \cup _{i=0}^{k}S_{i}
\]
where $\mcS(\bar{\bon})=(S_0,\ldots)$ is the family of indices 
associated to the sequence $(\mathbf{n}_{k})$.
\end{lemma}

\begin{proof}
We inductively find indices $\bo{n}_{1}$, $\bo{n}_{2},\ldots$,
as above, open sets $V\supseteq V_{1}\supseteq V_{2} \supseteq\ldots$
and points $x_1, x_2, \ldots$ with $x_i \in V_i$ and
$x_i=\phi_{\bo{k}_i}(x)$ for some index $\bo{k}_i$, so that
\[
\phi_{\bo{s}}(V_{i})\subseteq V\quad \text{for all}
    \quad \bo{s}\in S_{i}.
\]
This will prove the Lemma, for if $k\in \mathbb{N}$ and $\mathbf{s}\in S_{i}$
for some $i\leq k$ then, since $x_{k}\in V_{k}\subseteq V_{i}$ it will
follow that $\phi _{\mathbf{s}}(x_{k})\in \phi _{\mathbf{s}}(V_{k})\subseteq
\phi _{\mathbf{s}}(V_{i})\subseteq V$.

\medskip
Since $\lim_{k}\phi_{\bo{p}_{k}}(x)=x\in V$, there is
$\bo{p}_{i_1}$ with $x_1=\phi_{\bo{p}_{i_1}}(x)\in V$. Let
$\bo{k}_1=\bo{p}_{i_1}$. Using $\lim_k\phi_{\bo{p}_k}(x)=x\in V$
and the form of the $\bo{p}_k$, it follows that there is $\bo{n}_1
\in \Delta_J$ so that $\phi_{\bo{n}_1+\bo{k}_1}(x)\in V$.
Now
\[ \phi_{\bo{n}_1}(x_{1}) =\phi_{\bo{n}_1}(\phi_{\bo{k}_1}(x))
   =\phi_{\bo{n}_1+\bo{k}_1}(x) \in V,
\]
and so there is $V_1\subseteq V$, an open neighbourhood of $x_1$,
so that $\phi_{\bo{n}_1}(V_1)\subseteq V$. Since
$S_1=\{\bo{n}_1\}$, this establishes the base step.

\medskip

For the inductive step, assume we have chosen indices
$\bo{n}_1,\bo{n}_2,\ldots,\bo{n}_q$, open subsets of $V$,
$V_1\supseteq V_2\supseteq\ldots\supseteq V_q$ and points $x_1,
x_2, \ldots, x_q$, with $x_i\in V_i$ and $x_i=\phi_{\bo{k}_i}(x)$,
so that, for $i=1,\ldots,q$, we have
\begin{equation} \label{i}
\phi_\bo{s}(V_i)\subseteq V\text{ for all }\bo{s}\in S_i.
\end{equation}

Since $\lim_{k}\phi_{\bo{p}_{k}}(x_{q}) = \phi_{\bo{k}_q}
\bigl(\lim_{k}\phi_{\bo{p}_{k}}(x)\bigr) = x_q\in V_q$, there is
$\bo{k}_{q+1}=\bo{p}_{i_q}$ so that
$x_{q+1}=\phi_{\bo{k}_{q+1}}(x_q)\in V_q$. Notice that $\bo{m}_q$
(as in Definition 6) is in $\Delta_J$. It follows that
there exists $\bo{n}_{q+1}\in \Delta_J$ such that
$\phi_{\bo{n}_{q+1}+\bo{m}_{q}+\bo{k}_{q+1}}(x_q)\in V_q$ and so
$\phi_{\bo{n}_{q+1}+\bo{m}_q}(x_{q+1})\in V_q$. Hence there exists
an open neighbourhood $V_{q+1}$ of $x_{q+1}$, contained in $V_q$,
so that
\begin{equation}
\phi_{\bo{n}_{q+1}+\bo{m}_q}(V_{q+1})\subseteq
V_{q}.\label{i1}%
\end{equation}

It remains only to show that $\phi_\bo{s}(V_{q+1})\subseteq V$ for
all $\bo{s}\in S_{q+1}$. An element $\bo{s}$ in $S_{q+1}$ is of
the form $\bo{s}=\bo{n}_{q+1}+\bo{m}_{q}+\bo{j}$ for some $\bo{j}
\in \cup_{i=0}^q S_i$. Assuming $\bo{j} \in S_{i}$ for some $i$,
we have
\begin{align*}
\phi_{\bo{s}}(V_{q+1})
    & =\phi_{\bo{j}}(\phi_{\bo{n}_{q+1}+\bo{m}_q}(V_{q+1})) \subseteq
    \phi_{\bo{j}}(V_q)\qquad &\text{by (\ref{i1})}\\
    & \subseteq\phi_{\bo{j}}(V_i)\subseteq V &\text{by (\ref{i})}
\end{align*}
completing the induction.
\end{proof}

\begin{proof}[Proof of Lemma~\ref{rad}]
Assume that $f(x) \ne0$ for some $J$-recurrent point $x$. We will
find $B \in\mcA$ such that $B U_{\bo{q}} f$ has nonzero spectral
radius. We may scale $f$ so that there exists a relatively compact
open neighbourhood $V$ of $x$ such that $|f(y)|\geq1$ for all
$y\in V$. Since $U_\bo{q} |f|^2 =(U_{\bo{q}}f)f^* \in\rad\mcA$
when $U_\bo{q} f\in\rad\mcA$, we may also assume that $f\geq0$.

Since $x$ is $J$-recurrent, there exists a sequence
$(\bo{p}_{k})$ which is strictly increasing in the directions
of $J$ such that $\lim_{k}\phi_{\bo{p}_{k}}(x)=x$.
Deleting some initial segment, we may assume that
$\phi_{\bo{p}_{k}}(x)\in V$ for all $k\in\mathbb{N}$.

If $(\bo{p}_k)$ has all entries going to infinity, then we may
apply Lemma \ref{rec} with $J=\{1,2,\ldots d\}$, to find a
strictly increasing sequence $(\bon_{k})$ such that
$\bon_k>\bo{q}$ for all $k$ and points $x_k \in V$ such that
$\phi_{\bo{s}}(x_k)\in V$ for all $\bo{s}$ in
$\cup_{i=0}^{k}S_{i}.$

If not, (enlarging $J$ and passing to a subsequence if necessary)
we may assume that the restriction of $(\bo{p}_{k})$ to $J^{c}$
takes only finitely many values. Passing to another subsequence,
we may further assume that this restriction is constant. Applying
Lemma~\ref{rec}, we may find a strictly increasing sequence
$(\bon_k)$ in $\bbZ_+^d$ with $\bon_k \in \Delta_J$ and points
$x_k \in V$ such that $\phi_{\bo{s}}(x_k)\in V$ for all $\bo{s}$
in $\cup_{i=0}^{k}S_{i}.$ We may suppose that
$\bon_{k}-\bo{q}\in\Sigma_{J}$ for all $k$. Thus
$U_{\bo{n}_k-\bo{q}}$ is an admissible term in the formal power
series of an element of $C_0(X)\times_\phi \bbZ_+^d$.

Fix a nonnegative function $h\in C_{0}(X)$ such that
$h(\phi_{\bo{q}}(y))=1$ for all $y\in V$ and consider
\[
B=\sum_{k=1}^{\infty} U_{\bon_{k}-\bo{q}} \frac{h}{2^{k-1}}.
\]
This is an element of $\mcA$ since the series converges
absolutely. To complete the proof, it suffices to show that the
spectral radius of $A\equiv BU_{\bo{q}}f$ is strictly positive.
Note that
\[ A =\sum U_{\bon_{k}} \frac{g}{2^{k-1}}, \]
where $g$ is $f.(h\circ\phi_{\bo{q}})$, a nonnegative function
satisfying $g(y)\geq1$ for all $y\in V$. Thus each Fourier
coefficient $E_{\bon}(A^{m})$ of $A^{m}$ is a finite sum of
nonnegative functions, and hence its norm dominates the (supremum)
norm of each summand. Since $\Vert A^{2^{k}-1}\Vert\geq\Vert
E_{\bon}(A^{2^{k}-1})\Vert$, it suffices to find $\varepsilon>0$
such that for each $k$ there exists $\bon$ such that the norm of
\emph{some} summand of $E_{\bon}(A^{2^{k}-1})$ exceeds
$\varepsilon^{2^{k}-1}$.

If we let $P_1=U_{\bon_1}g$, then trivially $P_1$ is a term
in $A$.
In the next product, $A^3=A(\sum U_{\bon_k} \frac{g}{2^{k-1}})A=
\sum A(U_{\bon_k} \frac {g}{2^{k-1}}) A$, we have the term
\[
P_2=U_{\bon_1} g \left( U_{\bon_2} \frac{g}{2} \right)
    U_{\bon_1} g.
\]
Generally, one term in the expansion of
$A^{2^k-1}=A^{2^{k-1}-1}AA^{2^{k-1}-1}$ is
\[
P_k=P_{k-1}\left(U_{\bon_k} \frac{g}{2^{k-1}}\right)P_{k-1}.
\]

\noindent\textbf{Claim.}
If $\lambda_1=1$ and $\lambda_{k+1}=\lambda_{k}^{2}/2^{k}$, then
$\displaystyle
P_{k}= U_{\bom_{k}} \lambda_{k} \prod_{\bo{s}}g\circ\phi_{\bo{s}}
$
where $\bom_{k}$ is as in the definition of $\mcS(\bar{\bon})$ and
the product is over all $\bo{s}$ in
$(\cup_{i=0}^{k}S_i)\backslash\{\bom_{k}\}$.

\medskip\noindent\emph{Proof of Claim}.
For $k=1$, the claim holds trivially
as $(S_0\cup S_1)\backslash\{\bom_1\}=\{0\}$. Assuming
the claim is true for some $k$, we have
\begin{align*}
P_{k+1} & =P_{k}\left(U_{\bon_{k+1}} \frac{g}{2^k}\right) P_k\\
        & =U_{\bom_{k}} \lambda_k \left( \prod_{\bo{s}}g\circ\phi_{\bo{s}}
        \right) \left(U_{\bon_{k+1}} \frac{g}{2^{k}}\right)
        U_{\bom_{k}} \lambda_{k}
        \left(\prod_{\bo{t}}g\circ\phi_{\bo{t}}\right) \\
\intertext{(where $\bo{s,t}$ range over
$(\cup_{i=0}^{k}S_i)\backslash\{\bom_{k}\}$) }
    &= U_{\bom_k} \frac{\lambda_k^2}{2^{k}}
        \left( \prod_{\bo{s}}g\circ\phi_{\bo{s}} \right)
        U_{\bon_{k+1}+\bom_k} \, (g\circ\phi_{\bom_k})
        \left( \prod_{\bo{t}}g\circ \phi_{\bo{s}} \right) \\
    &= U_{2\bom_{k}+\bon_{k+1}} \frac{\lambda_{k}^{2}}{2^{k}}
        \left( \prod_{\bo{s}}g\circ\phi_{\bo{s}+\bon_{k+1}+\bom_k}\right)
        (g\circ\phi_{\bom_{k}})
        \left( \prod_{\bo{t}}g\circ\phi_\bo{t} \right) \\
    & = U_{2\bom_{k}+\bon_{k+1}} \frac{\lambda_{k}^{2}}{2^{k}}
        \left( \prod_{\bo{s}^{\prime}}g\circ \phi_{\bo{s}^{\prime}}\right)
        \left( \prod_{\bo{t}^{\prime}}g\circ \phi_{\bo{t}^{\prime}}\right)
\intertext{where $\bo{s}^{\prime}$ ranges over $\{\bon_{k+1}+\bom_{k}+\bo{s}\}$,
for $\bo{s}\in(\cup_{i=0}^{k}S_i)\backslash \{\bom_{k}\}$, and
$\bo{t}^{\prime}$ ranges over $(\cup_{i=0}^{k}S_i)$.
Therefore}
P_{k+1} &= U_{\bom_{k+1}} \lambda_{k+1}
    \left( \prod_{\bo{s}}g\circ\phi_{\bo{s}}\right)
\end{align*}
for $\bo{s}$ in $(\cup_{i=0}^{k+1}S_i)\backslash\{\bom_{k+1}\}$,
proving the claim.

\medskip

Recall that for each $k\in \mathbb{N}$ there exists $x_{k}\in V$
such that $\phi _{\bo{s}}(x_{k})\in V$ for all
$\bo{s}\in \cup _{i=0}^{k}S_{i}$. Since $g|_{V}\geq 1$, we have
$\displaystyle\prod_{\bo{s}}g(\phi _{\bo{s}}(x_{k}))\geq 1$ where
$\bo{s}$ ranges over  $(\cup_{i=0}^{k}S_i)\backslash \{\bo{m}_k \}$
and hence
$\displaystyle\Vert\prod_{\bo{s}}g\circ \phi _{\bo{s}}\Vert \geq 1$.
 From the claim, it follows that $\Vert P_k\Vert \geq \lambda _k$
and so, by the earlier remarks,
\[
\Vert A^{2^{k}-1}\Vert\geq\Vert
E_{\bom_{k}}(A^{2^{k}-1})\Vert\geq\Vert
P_{k}\Vert\geq\lambda_{k}.
\]
Thus the proof will be complete if we show that $\lambda_{k}\geq(\frac{1}%
{2})^{2^{k}-1}$ or equivalently
$\log_2\lambda_{k}^{-1}\leq2^{k}-1$ for all $k$.
Setting $\mu_{k}=\log_2\lambda_{k}^{-1}$, the recurrence
relation for $\lambda_k$ becomes $\mu _{k+1}=2\mu_{k}+k$ and $\mu_1=0$,
which has solution $\mu_{k}=2^{k}-k-1$.
\end{proof}

\section{Wandering sets and Semisimplicity}

We characterise semisimplicity of analytic crossed products and
show this is equivalent to being semiprime.
Part of this characterisation is of course a special case of our
main result, Theorem~\ref{fulrad2}, but we will need the
preliminary results in any case.

A \emph{wandering open set} is an open set $V\subset X$ so that
$\phi_{\bon}^{-1}(V)\cap V= \emptyset$ whenever $\bo{n}\in\bbZ_+^d$ is nonzero.
A \emph{wandering point }is a point with a wandering neighbourhood.

We will need the following variant: let $J \subseteq \{ 1,\ldots d
\}$. An open set $V\subseteq X$ is said to be \emph{wandering in
the directions of $J$}, or \emph{$J$-wandering}, if
$\phi_{\bon}^{-1}(V)\cap V=\emptyset$ whenever $\bo{n}$ is in $\Sigma_J$.
It is easily seen that, if $X_{Jw}$ denotes the set of all $J$-wandering
points (those with a $J$-wandering neighbourhood), then $X_{Jw}$ is open
and its complement is invariant and contains the set $X_{Jr}$ of
$J$-recurrent points.

Note, however, that it is possible for a recurrent point (in the
usual sense) to have a neighbourhood which is $J$-wandering (for
some $J$). For example, if $X = \mathbb{R}^2$ and
$\phi_1(x,y)=(x+1,y)$ while $\phi_2(x,y)=(x,3y)$, then the origin
is recurrent for the dynamical system $(X, (\phi_1,\phi_2))$, but
it also has a $\{1\}$-wandering neighbourhood.

\medskip

The idea of the following Lemma comes from~\cite[Theorem~4.2]{Muh83}.

\begin{lemma}\label{prad}
Suppose $V\subseteq X$ is an open set which is $J$-wandering and
$g\in C_{0}(X)$ is a nonzero function with support contained in
$V$. If $\bo{e}_J$ denotes the characteristic function of $J$,
then $B=U_{\bo{e}_J} g$ generates a nonzero ideal $\mcA B\mcA$
whose square is $0$.
\end{lemma}

\begin{proof} Let $C\in\mathcal{A}$ be arbitrary and
$h = E_{\bo{k}}(C)$.
Then
\[ BU_{\bo{k}} E_{\bo{k}}(C)B
  = U_{\bo{e}_J} g U_\bo{k} h U_{\bo{e}_J}g
  = U_{\bo{k+}2\bo{e}_J}
    (\alpha_{\bo{k+e}_J}(g) \alpha_{\bo{e}_J}(h)g)
\]
which is zero since $g$ is supported on $V$ and
$\alpha_{\bo{k+e}_J}(g)$ is supported on the disjoint set
$\phi_{\bo{k+e}_J}^{-1}(V)$.
This shows that all Fourier coefficients of $BCB$ will vanish,
and hence $BCB=0$.
It follows that all products $(C_1BC_2)(C_3BC_{4})$ vanish and hence
$(\mathcal{A}B\mathcal{A})^{2}=0$. On the other hand, choosing
functions $h_1\in C_{0}(X)$ equal to $1$ on
$\phi_{\bo{e}_J}^{-1}(V)$ and $h_2$ equal to $1$ on $V$, we find
$E_{\bo{e}_J}(h_1Bh_2)=\alpha_{\bo{e}_J}(h_1) g h_2=g\neq 0$, so
the ideal $\mathcal{A}B\mathcal{A}$ is nonzero.
\end{proof}

The following proposition is known for the usual notions of
recurrence and wandering in the case $d=1$;
see~\cite[Theorem~1.27]{Furstenberg81}.

\begin{proposition}
\label{wandrec} Suppose $X$ is a locally compact metrisable space.
If $(X,\Phi)$ has no nonempty $J$-wandering open
sets, then the $J$-recurrent points are dense.
\end{proposition}

\begin{proof}
Let $V\subseteq X$ be a relatively compact open set. We wish to
find a $J$-recurrent point in $V$.

Since $V$ is not $J$-wandering, there exists $\bon_1\in\Sigma_J$
such that $\phi_{\bon_1}^{-1}(V)\cap V\neq\emptyset$.
Hence there is a nonempty, relatively compact, open set $V_1$
with $\mathrm{diam}(V_1)<1$ such that
$\overline{V_1}\subseteq\phi _{\bon_1}^{-1}(V)\cap V$.

Since $V_1$ contains no $J$-wandering subsets, a similar argument
shows that there exists $\bon_2$ such that
$\phi_{\bon_2}^{-1}(V_1)\cap V_1\neq\emptyset$ and the $j$-th
entry of $\bon_2$ is greater than that of $\bon_1$ for every $j\in
J$.

Inductively one obtains a sequence of open sets $V_{k}$ and
$\bon_{k}$ strictly increasing in the directions of $J$ with
$\overline{V_{k}}\subseteq
\phi_{\bon_{k}}^{-1}(V_{k-1})\cap V_{k-1}$ and $\mathrm{diam}%
(V_{k})<1/k$ all contained in the compact metrisable space
$\overline{V_{0}}$. It follows from Cantor's theorem that the
intersection $\cap_{n\geq1}\overline{V_{n}}$ is a singleton, say
$x$. Since
$x\in\overline{V_{k}}\subseteq\phi_{\bon_{k}}^{-1}(V_{k-1})$ we
have $\phi_{\bon_{k}}(x)\in V_{k-1}$ for all $k$ and so $\phi
_{\bon_{k}}(x)\rightarrow x$; hence $x\in X_{Jr}$.
\end{proof}

\begin{theorem} \label{sem}
If $X$ is a metrisable, locally compact space, then the following
are equivalent:
\begin{enumerate}
\item the strongly recurrent points are dense in $X$,
\item $C_0(X)\times_\phi\bbZ_+^d$ is semisimple, and
\item $C_0(X)\times_\phi\bbZ_+^d$ is semiprime.
\end{enumerate}
\end{theorem}

\begin{proof}
If the strongly recurrent points are dense in $X$, then by
Lemma~\ref{rad} there are no nonzero monomials in the Jacobson
radical of $C_0(X)\times_\phi\gS$. But we have already observed
that an element  $A$ is in the Jacobson radical if and only if
each monomial $U_\bon E_\bon(A)$ is. Thus $C_0(X)\times_\phi\gS$
is semisimple and hence semiprime.

Suppose that $C_0(X)\times_\phi\gS$ is semiprime. Then Lemma
\ref{prad} shows that there are no nonempty $J$-wandering open
sets for $J= \{1,2, \ldots,d \}$
Thus, by Proposition~\ref{wandrec}, the strongly recurrent points
are dense.
\end{proof}

\section{Centres and the Jacobson Radical}

In order to describe the Jacobson radical of an analytic crossed
product, we need to characterise the closure of the $J$-recurrent
points, for a dynamical system $(X,\Phi)$ with $X$ a locally
compact metrisable space.

\begin{lemma}\label{cent1}
$(i)$ If $Y\subseteq X$ is a closed invariant set, the set
$Y_{Jr}$ of $J$-recurrent points  for the dynamical system
$(Y,\Phi)$ equals $X_{Jr}\cap Y$.

$(ii)$ The set $\overline{X_{Jr}}$ is the largest closed invariant
set $Y\subseteq X$ such that $(Y,\Phi)$ has no $J$-wandering
points.
\end{lemma}

\begin{proof}
(i) To see
that $Y_{Jr}\subseteq X_{Jr}$, note that if $y\in Y_{Jr}$ then for
every neighbourhood $V$ of $y$ (in $X$) the set $V\cap Y$ is a
neighbourhood of $y$ in the relative topology of $Y$, so there
exists $\bon\in\gS_{J}$ such that $\phi_{\bon}(y)\in V\cap Y$.
Thus $\phi_{\bon}(y)\in V$ showing that $y\in X_{Jr}$. On the
other hand if $y\in Y\cap X_{Jr}$ then for each relative
neighbourhood $V\cap Y$ of $y$, since $V$ is a neighbourhood of
$y$ in $X$ there exists $\bon\in\gS_{J}$ such that
$\phi_{\bon}(y)\in V$. Since $y\in Y$ and $Y$ is invariant,
$\phi_{\bon}(y)\in V\cap Y$ establishing (i).

(ii) Given a closed invariant set $Y \subseteq X$, if $(Y,\Phi|_Y)$ has
no $J$-wandering points, then $Y_{Jr}$ is dense in $Y$ by
Proposition~\ref{wandrec}, and hence
$Y\subseteq\overline{X_{Jr}}$. On the other hand,
$(\overline{X_{Jr}},\Phi)$ clearly has no $J$-wandering open
sets.
\end{proof}

The set $\overline{X_{Jr}}$ is found by successively `peeling off'
the $J$-wandering parts of the dynamical system. This construction
and Lemma~\ref{cent2} generalises the well
known concept of the centre of a dynamical system
$(X,\phi)$~\cite[7.19]{GottschalkHed55}.

If $V\subseteq X$ is the union of the $J$-wandering open subsets
of $X$, then let $X_{J,1}$ be the closed invariant set
$X\backslash V$. Consider the dynamical system $(X_{J,1},\Phi_{J,1})$,
where $\Phi_{J,1}\equiv\Phi|_{X_{J,1}}$. Let $X_{J,2}$ be the
complement of the union of all $J$-wandering open sets of
$(X_{J,1},\Phi_{J,1})$. Again we have a closed invariant set, and we
may form the dynamical subsystem $(X_{J,2},\Phi_{J,2})$ where
$\Phi_{J,2}\equiv \Phi|_{X_{J,2}}$. By transfinite recursion, we obtain a
decreasing family $(X_{J,\gamma},\Phi_{J,\gamma})$ of dynamical
systems: indeed, if $(X_{J,\gamma},\Phi_{J,\gamma})$ has been
defined, we let $X_{J,\gamma +1}\subseteq X_{J,\gamma}$ be the set
of points in $(X_{J,\gamma},\Phi_{J,\gamma})$ having no
$J$-wandering neighbourhood and we define
$\Phi_{J,\gamma+1}=\Phi|_{X_{J,\gamma+1}}$; if $\beta$ is a limit
ordinal and the systems $(X_{J,\gamma},\Phi_{J,\gamma})$ have been
defined for all $\gamma<\beta$ then we set
$X_{J,\beta}=\cap_{\gamma<\beta }X_{J,\gamma}$ and
$\Phi_{J,\beta}=\Phi|_{X_{J,\beta}}$. (We write $X_{J,0}=X$ and
$\Phi_{J,0}=\Phi$.) This process must stop, for the cardinality of
the family $\{X_{J,\gamma}\}$ cannot exceed that of the power set
of $X$.

\begin{definition}
By the above argument, there exists a least ordinal $\gamma$
such that $X_{J,\gamma+1}=X_{J,\gamma}$.
The set $X_{J,\gamma}$ is called \emph{the strong
$J$-centre} of the dynamical system, and $\gamma$ is called
\emph{the depth} of the strong $J$-centre.
\end{definition}

\begin{lemma} \label{cent2}
If $X$ is metrisable, then the strong $J$-centre of the dynamical
system is the closure of the $J$-recurrent points.
\end{lemma}

\begin{proof}
As a $J$-recurrent point cannot be $J$-wandering,
$X_{Jr}\subseteq X_{J,1}$. If $X_{Jr}\subseteq X_{J,\gamma}$ for some
$\gamma$, then by Lemma \ref{cent1} the set $(X_{J,\gamma})_{Jr}$ of
$J$-recurrent points of the subsystem
$(X_{J,\gamma},\Phi_{J,\gamma})$ equals $X_{Jr}\cap X_{J,\gamma}$, so
$(X_{J,\gamma})_{Jr}=X_{Jr}$; but $(X_{J,\gamma})_{Jr}\subseteq
X_{J,\gamma+1}$, and so $X_{Jr}\subseteq X_{J,\gamma+1}$. Finally, if
$\gamma$ is a limit ordinal and we assume that $X_{Jr}\subseteq
X_{J,\delta}$ for all $\delta<\gamma$ then
$X_{Jr}\subseteq\cap_{\delta<\gamma}X_{J,\delta}=X_{J,\gamma}$. This
shows that $X_{Jr}\subseteq\cap_{\gamma}X_{J,\gamma}$ and so
$\overline{X_{Jr}}\subseteq \cap_{\gamma}X_{J,\gamma}$ since the sets
$X_{J,\gamma}$ are closed.

But on the other hand, if $\gamma_0$ is the depth of the strong
$J$-centre we have $\cap_{\gamma}X_{J,\gamma}=X_{J,\gamma_0}$, a closed
invariant set. Since $X_{J,\gamma_0+1}=X_{J,\gamma_0}$, the dynamical
system $(X_{J,\gamma_0},\Phi_{J,\gamma_0})$ can have no $J$-wandering points.
Thus it follows from Lemma \ref{cent1} that
$X_{J,\gamma_0}\subseteq\overline{X_{Jr}}$ and hence equality holds.
\end{proof}

\begin{remark}
If $X$ is a locally compact (not necessarily metrisable) space
and $\{ \phi_{\bon}: \bon \in \bbZ^d \}$ is an action of an
\emph{equicontinuous} group of homeomorphisms (with respect to a uniformity
compatible with the topology of $X$) then $X_{Jr}=X\setminus X_{Jw}$
(see \cite{amanou}, Proposition 4.15).
\end{remark}

\begin{lemma}
\label{prad4} For any ordinal $\delta ,$ any $f\in
C_{c}(X_{J,\delta +1}^{c})$ $($i.e. $f$ has compact support
disjoint from $X_{J,\delta +1})$ can be written as a \emph{finite}
sum $f=\sum f_{k}$ where each $f_{k}$ has compact support
contained in a set $V_k$ such that $V_k \cap
X_{J,\delta }$ is $J$-wandering set for $(X_{J,\delta },\Phi
_{J,\delta })$.
\end{lemma}

\begin{proof}
If $K$ is the support of $f$ then $K\cap X_{J,\delta}\subseteq
X_{J,\delta }\setminus X_{J,\delta +1}$; in other words the
compact set $K\cap X_{J,\delta }$ consists of $J$-wandering
points for $(X_{J,\delta },\Phi _{J,\delta })$. This means that each
$x\in K\cap X_{J,\delta }$ has an open neighbourhood $V_{x}$ so that
the (relatively open) set $V_{x}\cap X_{J,\delta }$ is
$J$-wandering for $(X_{J,\delta },\Phi _{J,\delta })$.
Each $y \in K\setminus X_{J,\delta}$ has an open neighbourhood
$V_y$ such that $V_y\cap X_{J,\delta}$ is empty (and so $J$-wandering).

The family $\{V_x:x\in K\}$ is an open cover for $K$.
Thus, there is a partition of unity for $f$, i.e., a finite subcover,
$\{V_k:1\leq k\leq m\}$, and functions $f_k$, $1 \le k \le m$,
with $\mathrm{supp} (f_k) $ a compact subset of $V_k$, so that
$f=f_1+\ldots +f_m$.
\end{proof}

\begin{definition}
We denote by $\mcR_{J,\gamma}$ the closed ideal generated by all
monomials of the form $U_\bon f$ where $\bon$ is in $\Sigma_J$ and
$f\in C_0(X)$ vanishes on the set $X_{J,\gamma}$ and by
$\mcS_{J,\gamma}$ the set of all elements of the form $B f$ where
$B\in \mcR_{J,\gamma}$ and $f$ has compact support disjoint from
$X_{J,\gamma}$. 
\end{definition}

Note that
a monomial $U_{\bon} f\in \mcR_{J,\gamma}$ may be written in the form
$CU_{\bo{e}_J}f$ with $C \in \mcA$, since $\bon \in \Sigma_J$.

Also observe that $\mcS_{J,\gamma }$ is dense in
$\mcR_{J,\gamma}$.
Indeed if $U_{\bon}f\in \mathcal{R}_{J,\gamma }$, then
$f$ can be approximated by some $g\in C_c(X_{J,\gamma}^c)$;
now $U_{\bon}g$ is in $\mcS_{J,\gamma}$ and approximates $U_{\bon}f$.

\begin{proposition}\label{indu}
For each ordinal $\gamma$ and each
$J\subseteq\{1,2,\ldots,d\}$, the set $\mcS_{J,\gamma}$
is contained in $\rad\mcA$.
Hence $\mathcal{R}_{J,\gamma }$ is contained in
$\rad\mcA$.

If $\prad\mcA$ is closed, then  $\mcR_{J,\gamma}$ is contained
in $\prad\mcA$.
\end{proposition}

\begin{proof}
Since $\mathcal{S}_{J,\gamma }$ is dense in $\mcR_{J,\gamma}$,
it suffices to prove that any $A=Bf\in \mathcal{S}_{J,\gamma} $
is contained in $\mathrm{Rad}\mathcal{A}$.

Suppose $\gamma =1.$ By Lemma \ref{prad4}  we may write
$A$ as a finite sum $A=\sum_k Bf_{k}$ where
each $f_{k}$ is supported on a compact set which is
$J$-wandering.
Since $A_{k}\equiv B f_{k}=D U_{\bo{e}_J} f_{k}$ for
some $D\in \mcA$ as observed above, by Lemma~\ref{prad} we have
$(\mathcal{A}A_{k}\mathcal{A})^2=0$ and so
$A_{k}\in \mathrm{PRad}\mathcal{A}$.
Thus $A\in \mathrm{PRad}\mathcal{A} \subseteq \rad\mcA$.

Suppose the result has been proved for all ordinals less than some
$\gamma$.

Let $\gamma $ be a limit ordinal.
If $\mathrm{supp}\, f=K\subseteq X_{J,\gamma }^{c}$, we have
$K\subseteq X_{J,\gamma }^{c}=
\cup_{\delta <\gamma}X_{J,\delta }^{c}$, hence $K$ can be covered by
finitely many of the $X_{J,\delta }^{c}$, hence (since they are
decreasing) by one of them. Thus $f$ has compact support contained
in some $X_{J,\delta }^{c}$ ($\delta <\gamma $) and so
$Bf \in \mathcal{S}_{J,\delta }$.
Therefore $A=Bf\in \mathrm{Rad}\mathcal{A}$ by the induction
hypothesis.

Now suppose that $\gamma$ is a successor, $\gamma =\delta +1$.
By Lemma~\ref{prad4}, we may write $f=\sum f_{k}$ where the support
of $f_{k}$ is compact and
contained in an open set $V_k$ such that
$V_{k}\cap X_{J,\delta }$ is $J$-wandering for
$(X_{J,\delta },\Phi_{J,\delta })$, i.e.,
\[ \phi _{\bon}^{-1}(V_{k}\cap X_{J,\delta })\cap
    (V_{k}\cap X_{J,\delta})=\emptyset
\]
when $\bon \in \Sigma_J$.
This can easily be seen to imply
$\phi_{\bon}^{-1}(V_{k})\cap
V_{k}\subseteq X_{J,\delta }^c$.

Let $C\in \mcA$ be arbitrary. Writing $A_k=D U_{\bo{e}_J} f_k$
as above, it follows as in the proof of Lemma~\ref{prad}
that for each $k$ all Fourier coefficients of $A_kCA_k$ are supported
in $V_k\cap \phi _{\bon}^{-1}(V_k)$ (for some $\bon \in\Sigma_J$)
which is contained in $X_{J,\delta }^{c}$ by the
previous paragraph.

Thus $A_{k}CA_{k}\in \mcR_{J,\delta}$.
By the induction hypothesis, $A_{k}CA_{k}$ must be contained in
$\rad\mcA$.
Thus $(A_{k}C)^2$ is quasinilpotent, hence so is
$A_{k}C$ (by the spectral mapping theorem). Since $C \in \mcA$ is
arbitrary, it follows that  $A_k \in \rad\mcA$ for each $k$, so
that $A\in \rad\mcA$.

Finally, we suppose that $\prad\mcA$ is closed.
Then the argument above can be repeated exactly up to the
previous paragraph, changing $\rad\mcA$ to $\prad\mcA$.
The previous paragraph can be replaced by the following
argument.

Thus $A_{k}CA_{k}\in \mathcal{R}_{J,\delta }$.
By the induction hypothesis, $A_{k}CA_{k}$ must be contained in
$\prad\mcA$. Thus all products $(C_1 A_{k}C_2)(C_3A_{k}C_{4})$ are
in $\prad\mcA$ and so the (possibly non-closed) ideal
$\mathcal{J}_{k}$ generated by $A_{k}$ satisfies
$\mathcal{J}_{k}\mathcal{J}_{k}\subseteq\prad\mcA$.
For every prime ideal $\mathcal{P}$, we have
$\mathcal{J}_{k}\mathcal{J}_{k}\subseteq \mathcal{P}$
and so $\mathcal{J}_{k}\subseteq \mathcal{P}$.
Hence $\mathcal{J}_{k}\subseteq \prad\mcA$, and therefore
$A_{k}\in \prad\mcA$ for each $k$, so that $A\in \prad\mcA$.
\end{proof}

One cannot conclude that $\mcR_{J,\gamma}\subseteq\prad\mcA$ in general,
even for finite $\gamma$, as the following example shows.
Thus the prime radical is not always closed.
Note that Hudson has given examples of TAF algebras in which
the prime radical is not closed~\cite[Example~4.9]{Hud97}.

\begin{example}
We use a continuous dynamical system $(X,\{\phi_{t}\}_{t \in
\bbR})$ based on~\cite[Example~3.3.4, p.~20]{BhSz} and look at the
discrete system given by the maps $\{\phi_t\}$ for $t \in \bbZ_+$.
The space $X$ is the closed unit disc in $\mathbb{R}^{2}$. For the
continuous system, the trajectories consist of: (i) three fixed
points, namely the origin $O$ and the points $A(1,0)$ and
$B(-1,0)$ on the unit circle, (ii) the two semicircles on the unit
circle joining $A$ and $B$ and (iii) spiraling trajectories
emanating at the origin and converging to the boundary.

Let $\phi =\phi_1$.
The recurrent points for the (discrete) dynamical system $(X,\phi)$ are
$X_{r}=\{A,B,O\}$ and the set of wandering points is the open unit disc
except the origin.
Hence $X_2=X_{r}$ and so the depth of the dynamical system is 2.

Now choose small disjoint open neighbourhoods $V_A,V_B,V_O$ around
the fixed points and let $f\in C(X)$ be a nonnegative function
which is $1$ outside these open sets and vanishes only at $A,B$ and $O$.
Then the element $Uf\in \mcA$ is clearly not nilpotent, so
$Uf\notin \prad\mcA$.
However $Uf\in \rad\mcA$ by the next theorem.
\end{example}

\begin{theorem} \label{fulrad2}
Let  $(X,\Phi)$ be a dynamical system with $X$ metrisable. The
Jacobson radical, $\rad(C_{0}(X)\times _{\phi }\bbZ_+^d)$, is
the closed ideal generated by all monomials $U_{\bon}f$ $(\bo{n
\ne 0})$ where $f$ vanishes on the set $X_{Jr}$ of $J$-recurrent
points corresponding to the support $J$ of $\bon$.

Moreover, $\prad\mcA = \rad\mcA$ if and only if $\prad\mcA$ is closed.
\end{theorem}

\begin{proof}
Let $U_{\bon}f$ be a monomial contained in $\rad\mcA$ and
let $J$ be the support of $\bon$. Then Lemma~\ref{rad} shows
that $f$ must vanish on $X_{Jr}$.

On the other hand, let $U_{\bon}f$ be as in the statement of the
Theorem, so that $f$ vanishes on $X_{Jr}$ (where
$J = \mathrm{supp}\,\bon$).
We will show that $U_{\bon}f$ is in $\rad\mcA$.
It is enough to suppose that the support $K$ of $f$ is compact.
Since $K$ is contained in
$(\overline{X_{Jr}})^c=\cup _{\gamma }X_{J,\gamma }^c$,
it is contained in finitely many, hence one,
$X_{J,\gamma }^{c}$. It follows by Proposition~\ref{indu} that
$U_{\bon}f\in \rad\mcA$.

In the final statement of the theorem, one direction is obvious.
For the other, suppose $\prad\mcA$ is closed.
Then by the final statement of Proposition~\ref{indu}, we have
$\mcR_{J,\gamma} \subseteq \prad\mcA$.
\end{proof}

This theorem leaves open the possibility that the closure of the
prime radical is always equal to the Jacobson radical.

\providecommand{\bysame}{\leavevmode\hbox
to3em{\hrulefill}\thinspace}

\end{document}